\documentclass[a4paper,10pt,fleqn]{article}
\usepackage[english]{babel}
\usepackage{amssymb}
\usepackage{graphicx}
\usepackage{natbib}
\usepackage{amsfonts}
\usepackage{amscd}
\usepackage{amsmath}
\usepackage{amsthm}
\usepackage{latexsym}
\usepackage[latin1]{inputenc}
\usepackage[T1]{fontenc}
\usepackage[pdftex,colorlinks,linkcolor=red,citecolor=blue]{hyperref}
\setcitestyle{notesep={; },authoryear,round,coma,aysep={,},yysep={;}}
\theoremstyle{plain}
\newtheorem{Thm}{Theorem}[section]
\newtheorem{Propo}[Thm]{Proposition}
\newtheorem{Lem}[Thm]{Lemma}
\newtheorem{Coro}[Thm]{Corollary}
\theoremstyle{definition}
\newtheorem{Def}{Definition}[section]
\theoremstyle{remark}
\newtheorem{Rem}{Remark}[section]
\newtheorem*{pf}{Proof}
\numberwithin{equation}{section}

\begin{document}
\title{Large deviations for backward stochastic differential equations driven by $G$-Brownian motion}
\author{Ibrahim DAKAOU\footnote{Email: idakaou@yahoo.fr} ~and Abdoulaye SOUMANA HIMA\footnote{Email: soumanahima@yahoo.fr}\\
\small Département de Mathématiques, Université Dan Dicko Dankoulodo de Maradi, \\
\small BP 465, Maradi, Niger}
\date{}
\maketitle
\begin{abstract}
In this paper, we consider forward-backward stochastic differential equation driven by $G$-Brownian motion ($G$-FBSDEs in short) with small parameter $\varepsilon > 0$. We study the asymptotic behavior of the solution of the backward equation and establish a large deviation principle for the corresponding process.


\end{abstract}

\section{Introduction}
\label{intro}
The large deviation principle (LDP in short) characterizes the limiting behavior, as $\varepsilon\rightarrow0$, of family of probability measures $\{\mu_{\varepsilon}\}_{\varepsilon>0}$ in terms of a rate function. Several authors have considered large deviations and obtained different types of applications mainly to mathematical physics. General references on large deviations are: \citet{Varadhan1984,Deuschel1989,Dembo1998}.

Let $X^{s,x,\varepsilon}$ be the diffusion process that is the unique solution of the following stochastic differential equation (SDE in short)
\begin{equation}\label{eq0}
X_{t}^{s,x,\varepsilon}=x+\int_{s}^{t}\beta(X_{r}^{s,x,\varepsilon})dr+\sqrt
\varepsilon \int_{s}^{t}\sigma(X_{r}^{s,x,\varepsilon})dW_{r},\; 0\leq s\leq t\leq T
\end{equation}
where $\beta$ is a Lipschitz function defined on $\mathbb{R}^{d}$ with values in $\mathbb{R}^{d}$, $\sigma$ is a Lipschitz function defined on $\mathbb{R}^{d}$ with values in $\mathbb{R}^{d\times k}$, and $W$ is a standard Brownian motion in $\mathbb{R}^{k}$ defined on a complete probability space $(\Omega,\mathcal{F},\mathbb{P})$. The existence and uniqueness of
the strong solution $X^{s,x,\varepsilon}$ of \eqref{eq0} is standard. Thanks to the work of \citet{Freidlin1984}, the sequence $(X^{s,x,\varepsilon})_{\varepsilon > 0}$ converges in probability, as $\varepsilon$ goes to 0, to $(\varphi_{t}^{s,x})_{s\leq t\leq T}$ solution  of the following deterministic equation
\begin{equation*}
\varphi_{t}^{s,x}=x+\int_{s}^{t}\beta(\varphi_{r}^{s,x})dr,\; 0\leq s\leq t\leq T
\end{equation*}
and satisfies a large deviation principle (LDP in short).

\citet{Rainero2006} extended this result to the case of backward stochastic differential equations (BSDEs in short) and \citet{Essaky2008} and \citet{Nzi2014a} to reflected BSDEs.

\citet{Gao2010} extended the work of \citet{Freidlin1984} to stochastic differential equations driven by $G$-Brownian motion ($G$-SDEs in short). The authors considered the following $G$-SDE: for every $0\leq t\leq T$,
\begin{equation*}
    X_{t}^{x,\varepsilon} =x+\int_{0}^{t}b^{\varepsilon}(X_{r}^{x,\varepsilon})dr+\varepsilon\int_{0}^{t}h^{\varepsilon}(X_{r}^{x,\varepsilon})d\langle B, B\rangle_{r/\varepsilon}+\varepsilon\int_{0}^{t}\sigma^{\varepsilon}(X_{r}^{x,\varepsilon})dB_{r/\varepsilon}
\end{equation*}
and use discrete time approximation to establish LDP for $G$-SDEs.

The aim of this paper is to establish LDP for $G$-BSDEs. More precisely, we consider the following forward-backward stochastic differential equation driven by $G$-Brownian motion: for every $s\leq t\leq T$,
\begin{equation*}
\begin{cases}
&X_{t}^{s,x,\varepsilon} =x+\int_{s}^{t}b(X_{r}^{s,x,\varepsilon})dr+\varepsilon\int_{s}^{t}h(X_{r}^{s,x,\varepsilon})d\langle B, B\rangle_{r}+\varepsilon\int_{s}^{t}\sigma(X_{r}^{s,x,\varepsilon})dB_{r} \\
&Y_{t}^{s,x,\varepsilon}=\Phi(X_{T}^{s,x,\varepsilon})+\int_{t}^{T}f(r,X_{r}^{s,x,\varepsilon},Y_{r}^{s,x,\varepsilon},Z_{r}^{s,x,\varepsilon})dr-\int_{t}^{T}Z_{r}^{s,x,\varepsilon}dB_{r}
\\ &\qquad\qquad+\int_{t}^{T}g(r,X_{r}^{s,x,\varepsilon},Y_{r}^{s,x,\varepsilon},Z_{r}^{s,x,\varepsilon})d\langle B, B\rangle_{r}-(K_{T}^{s,x,\varepsilon}-K_{t}^{s,x,\varepsilon})
\end{cases}
\end{equation*}
We study the asymptotic behavior of the solution of the backward equation and establish a LDP for the corresponding process.

The remaining part of the paper is organized as follows. In Section~\ref{Sec:Pre}, we present some preliminaries that are useful in this paper. Section~\ref{Sec:G-SDEs} is devoted to the large deviations for stochastic differential equations driven by $G$-Brownian motion obtained by \citet{Gao2010}. The large deviations for backward stochastic differential equations driven by $G$-Brownian motion are given in Section~\ref{Sec:G-BSDEs}.

\section{Preliminaries}
\label{Sec:Pre}
We review some basic notions and results about $G$-expectation, $G$-Brownian motion and $G$-stochastic integrals \citep[see][for more details]{Peng2010,Hu2014a}.

Let $\Omega$ be a complete separable metric space, and let $\mathcal{H}$ be a linear space of real-valued functions defined on $\Omega$ satisfying: if $X_i\in\mathcal{H}$, $i=1, \ldots, n$, then
\begin{equation*}
\varphi(X_1, \ldots, X_n)\in\mathcal{H}, \quad \forall\varphi\in\mathcal{C}_{l, Lip}(\mathbb{R}^n),
\end{equation*}
where $\mathcal{C}_{l,Lip}(\mathbb{R}^{n})$ is the space of real continuous functions defined on $\mathbb{R}^{n}$ such that for some $C>0$ and $k\in\mathbb{N}$ depending on $\varphi$,
\begin{equation*}
    \vert\varphi (x)-\varphi (y)\vert\leq C(1+\vert x\vert^{k}+\vert y\vert^{k})\vert x-y\vert,\quad \forall x,y\in\mathbb{R}^{n}.
\end{equation*}
\begin{Def}(\emph{Sublinear expectation space}).
A sublinear expectation $\widehat{\mathbb{E}}$ on $\mathcal{H}$ is a functional $\widehat{\mathbb{E}}:\mathcal{H}\longrightarrow \mathbb{R}$ satisfying the following properties: for all $X, Y\in\mathcal{H}$, we have
\begin{enumerate}
  \item Monotonicity: if $X\geq Y$, then $\widehat{\mathbb{E}}[X]\geq \widehat{\mathbb{E}}[Y]$;
  \item Constant preservation: $\widehat{\mathbb{E}}[c]=c$;
  \item Sub-additivity: $\widehat{\mathbb{E}}[X+Y]\leq \widehat{\mathbb{E}}[X]+\widehat{\mathbb{E}}[Y]$;
  \item Positive homogeneity: $\widehat{\mathbb{E}}[\lambda X]=\lambda\widehat{\mathbb{E}}[X]$, for all $\lambda\geq 0$.
\end{enumerate}
$(\Omega, \mathcal{H}, \widehat{\mathbb{E}})$ is called a \emph{sublinear expectation space}.
\end{Def}
\begin{Def}(\emph{Independence}).
Fix the sublinear expectation space $(\Omega ,\mathcal{H},\widehat{\mathbb{E}})$. A random variable $Y\in \mathcal{H}$ is said to be independent of $(X_{1}, X_{2}, \ldots, X_{n})$, $X_{i}\in\mathcal{H}$, if for all $\varphi\in\mathcal{C}_{l, Lip}(\mathbb{R}^{n+1})$,
\begin{equation*}
    \widehat{\mathbb{E}}\left[\varphi (X_{1}, X_{2}, \ldots, X_{n}, Y)\right]=
    \widehat{\mathbb{E}}\left[\widehat{\mathbb{E}}\left[\varphi(x_{1}, x_{2}, \ldots, x_{n}, Y)\right]\big|_{(x_{1}, x_{2}, \ldots, x_{n})=(X_{1}, X_{2}, \ldots, X_{n})}\right].
\end{equation*}
\end{Def}

Now we introduce the definition of $G$-normal distribution.
\begin{Def}(\emph{$G$-normal distribution}).
A random variable $X\in\mathcal{H}$ is called $G$-normally distributed, noted by $X\sim\mathcal{N}(0, [\underline{\sigma}^{2}, \overline{\sigma}^{2}])$, $0\leq\underline{\sigma}^{2}\leq\overline{\sigma}^{2}$, if for any function $\varphi\in\mathcal{C}_{l, Lip}(\mathbb{R})$, the fonction $u$ defined by $u(t,x):=\widehat{\mathbb{E}}[\varphi(x+\sqrt{t}X)],\;\left(t, x\right)\in\left[0,\infty\right)\times\mathbb{R}$, is a viscosity solution of the following $G$-heat equation:
\begin{equation*}\label{heat}
    \partial_{t}u-G\left(D_{x}^{2}u\right)=0,\;u(0,x)=\varphi(x),
\end{equation*}
where
\begin{equation*}
    G(a):=\frac{1}{2}(\overline{\sigma}^{2}a^{+}-\underline{\sigma}^{2}a^{-}).
\end{equation*}
In multi-dimensional case, the function $G(\cdot)$: $\mathbb{S}_{d}\longrightarrow\mathbb{R}$  is defined by
\begin{equation*}
    G(A)=\frac{1}{2}\sup_{\gamma\in\Gamma}\,\textrm{tr}(\gamma\gamma^{\tau}A),
\end{equation*}
where $\mathbb{S}_{d}$ denotes the space of $d\times d$ symmetric matrices and $\Gamma$ is a given nonempty, bounded and closed subset of $\mathbb{R}^{d\times d}$ which is the space of all $d\times d$ matrices.
\end{Def}
Throughout this paper, we consider only the non-degenerate case, i.e., $\underline{\sigma}^{2}>0$.

Let $\Omega:=\mathcal{C}\left([0, \infty)\right)$, which equipped with the raw filtration $\mathcal{F}$ generated by the canonical process $(B_t)_{t\geq 0}$, i.e., $B_t(\omega)=\omega_t$, for $(t, \omega)\in[0, \infty)\times\Omega$. Let us consider the function spaces defined by
\begin{align*}
  Lip(\Omega_{T})&:=\Big\{\varphi(B_{t_{1}}, B_{t_{2}}-B_{t_{1}}, \ldots, B_{t_{n}}-B_{t_{n-1}}): n\geq 1,\\
  &\quad\quad\quad\quad 0\leq t_{1}\leq t_{2}\leq\ldots\leq t_{n}\leq T, \varphi\in\mathcal{C}_{l,Lip}(\mathbb{R}^{n})\Big\},\quad \text{for}\quad T>0, \\
  Lip(\Omega)&:=\bigcup_{n=1}^{\infty}Lip(\Omega_{n}).
\end{align*}
\begin{Def}(\emph{$G$-Brownian motion and $G$-expectation}).
On the sublinear expectation space $\left(\Omega, Lip(\Omega), \widehat{\mathbb{E}}\right)$, the canonical process $(B_t)_{t\geq 0}$ is called a $G$-Brownian motion if the following properties are verified:
\begin{enumerate}
  \item $B_{0}=0$
  \item For each $t, s\geq 0$, the increment $B_{t+s}-B_t\sim\mathcal{N}(0, [s\underline{\sigma}^{2}, s\overline{\sigma}^{2}])$ and is independent from $(B_{t_{1}}, \ldots, B_{t_{n}})$, for $0\leq t_1\leq \ldots\leq t_n\leq t$.
\end{enumerate}
Moreover, the sublinear expectation $\widehat{\mathbb{E}}$ is called \emph{$G$-expectation}.
\end{Def}
\begin{Rem}\label{l6}
For each $\lambda>0$, $\left(\sqrt{\lambda}B_{t/\lambda}\right)_{t\geq 0}$ is also a $G$-Brownian motion. This is the \emph{scaling property} of $G$-Brownian motion, which is the same as that of the classical Brownian motion.
\end{Rem}
\begin{Def}(\emph{Conditional $G$-expectation}).
For the random variable $\xi\in Lip(\Omega_T)$ of the following form:
\begin{equation*}
\varphi(B_{t_{1}}, B_{t_{2}}-B_{t_{1}}, \ldots, B_{t_{n}}-B_{t_{n-1}}), \quad \varphi\in\mathcal{C}_{l, Lip}(\mathbb{R}^n),
\end{equation*}
the conditional $G$-expectation $\widehat{\mathbb{E}}_{t_{i}}[\cdot]$, $i=1, \ldots, n$, is defined as follows
\begin{equation*}
    \widehat{\mathbb{E}}_{t_{i}}[\varphi(B_{t_{1}}, B_{t_{2}}-B_{t_{1}}, \ldots, B_{t_{n}}-B_{t_{n-1}})]=\widetilde{\varphi}(B_{t_{1}},
    B_{t_{2}}-B_{t_{1}}, \ldots, B_{t_{i}}-B_{t_{i-1}}),
\end{equation*}
where
\begin{equation*}
    \widetilde{\varphi}\left(x_{1}, \ldots, x_{i}\right)=\widehat{\mathbb{E}}\left[\varphi\left(x_{1}, \ldots, x_{i}, B_{t_{i+1}}-B_{t_{i}}, \ldots, B_{t_{n}}-B_{t_{n-1}}\right)\right].
\end{equation*}
If $t\in(t_{i}, t_{i+1})$, then the conditional $G$-expectation $\widehat{\mathbb{E}}_{t}[\xi]$ could be defined by reformulating $\xi$ as
\begin{equation*}
    \xi=\widehat{\varphi}(B_{t_{1}}, B_{t_{2}}-B_{t_{1}}, \ldots, B_{t}-B_{t_{i}}, B_{t_{i+1}}-B_{t}, \ldots, B_{t_{n}}-B_{t_{n-1}}),
    \quad\widehat{\varphi}\in\mathcal{C}_{l,Lip}(\mathbb{R}^{n+1}).
\end{equation*}
\end{Def}

For $\xi\in Lip(\Omega_T)$ and $p\geq 1$, we consider the norm $\Vert\xi\Vert_{L^p_G}:=\left(\widehat{\mathbb{E}}\Big[\vert\xi\vert^p\Big]\right)^{1/p}$. Denote by $L^p_G(\Omega_T)$ the Banach completion of $Lip(\Omega_T)$ under $\Vert\cdot\Vert_{L^p_G}$. It is easy to check that the conditional $G$-expectation $\widehat{\mathbb{E}}_{t}[\cdot]: Lip(\Omega_T)\longrightarrow Lip(\Omega_t)$ is a continuous mapping and thus can be extended to $\widehat{\mathbb{E}}_{t}[\cdot]: L^p_G(\Omega_T)\longrightarrow L^p_G(\Omega_t)$.
\begin{Def}(\emph{$G$-martingale}).
A process $M=\left(M_{t}\right)_{t\in\left[0, T\right]}$ with $M_{t}\in L_{G}^{1}(\Omega_{t})$, $0\leq t\leq T$, is called a $G$-martingale if for all $0\leq s\leq t\leq T$, we have
\begin{equation*}
    \widehat{\mathbb{E}}_{s}[M_{t}]=M_{s}.
\end{equation*}
The process $M=\left(M_{t}\right)_{t\in\left[0, T\right]}$ is called symmetric $G$-martingale if $-M$ is also a $G$-martingale.
\end{Def}
\begin{Thm}\citep[Representation theorem of $G$-expectation, see][]{Hu2009,Denis2011}.
There exists a weakly compact set $\mathcal{P}\subset\mathcal{M}_1(\Omega_T)$, the set of probability measures on $(\Omega_T, \mathcal{B}(\Omega_T))$, such that
\begin{equation*}
    \widehat{\mathbb{E}}[\xi]:=\sup_{P\in \mathcal{P}}E_{P}[\xi] \quad \text{for all} \quad \xi\in L^1_G(\Omega_T).
\end{equation*}
$\mathcal{P}$ is called a set that represents $\widehat{\mathbb{E}}$.
\end{Thm}
Let $\mathcal{P}$ be a weakly compact set that represents $\widehat{\mathbb{E}}$. For this $\mathcal{P}$, we define the \emph{capacity} of a measurable set $A$ by
\begin{equation*}
    \widehat{C}(A):=\sup_{P\in\mathcal{P}}P(A), \quad A\in\mathcal{B}(\Omega_{T}).
\end{equation*}
A set $A\in\mathcal{B}(\Omega_{T})$ is a polar if $\widehat{C}(A)=0$. A property holds \emph{quasi-surely} (q.s.) if it is true outside a polar set.

An important feature of the $G$-expectation framework is that the quadratic variation $\left\langle B\right\rangle$ of the $G$-Brownian motion is no longer a deterministic process, which is given by
\begin{equation*}
    \left\langle B\right\rangle_{t}:=\lim_{\delta\left(\pi_{t}^{N}\right)
    \rightarrow 0}\sum_{j=0}^{N-1}(B_{t_{j+1}^{N}}-B_{t_{j}^{N}})^{2},
\end{equation*}
where $\pi_{t}^{N}=\{t_{0}, t_{1}, \ldots, t_{N}\}$, $N=1, 2, \ldots$, are refining partitions of $[0, t]$. By \cite{Peng2010}, for all $s, t\geq 0$, $\langle B\rangle_{t+s}-\langle B\rangle_{t}\in[s\underline{\sigma}^2, s\overline{\sigma}^2]$, $q.s.$

Let $M_{G}^{0}\left(0, T\right)$ be the collection of processes in the following form: for a given partition $\pi_{T}^{N}:=\{t_{0}, t_{1}, \ldots, t_{N}\}$ of $[0, T]$,
\begin{equation}\label{simp}
\eta_{t}\left(\omega\right)=\sum_{j=0}^{N-1}\xi_{j}\left(\omega\right)
\mathbf{1}_{\left[t_{j}, t_{j+1}\right)}(t),
\end{equation}
where $\xi_{i}\in Lip(\Omega_{t_{i}})$, for all $i=0, 1, \ldots, N-1$. For $p\geq 1$ and $\eta\in M_{G}^{0}\left(0, T\right)$, let $\left\Vert\eta\right\Vert_{H_{G}^{p}}:=\left(\widehat{\mathbb{E}}\left[\left(\int_{0}^{T}|\eta_{s}|^{2}ds\right)^{p/2}\right]    \right)^{1/p}$, $\Vert\eta\Vert_{M_{G}^{p}}:=\left(\widehat{\mathbb{E}}\left[\int_{0}^{T}|\eta_{s}|^{p}ds\right]\right)^{1/p}$ and denote by $H_{G}^{p}\left(0,T\right)$, $M_{G}^{p}(0,T)$ the completions of $M_{G}^{0}\left(0,T\right)$ under the norms $\Vert\cdot\Vert_{H_{G}^{p}}$, $\Vert\cdot\Vert_{M_{G}^{p}}$ respectively.

Let $\mathcal{S}_{G}^{0}\left(0,T\right):=\{h(t, B_{t_{1}\wedge t},B_{t_{2}\wedge t}-B_{t_{1}\wedge t}, \ldots, B_{t_{n}\wedge
t}-B_{t_{n-1}\wedge t}): 0\leq t_{1}\leq t_{2}\leq \ldots\leq t_{n}\leq T, ~h\in\mathcal{C}_{b,Lip}(\mathbb{R}^{n+1})\}$, where $\mathcal{C}_{b,Lip}(\mathbb{R}^{n+1})$ is the collection of all bounded and Lipschitz functions on $\mathbb{R}^{n+1}$. For $p\geq 1$ and $\eta\in\mathcal{S}_{G}^{0}\left(0, T\right)$, we set $\left\Vert\eta\right\Vert_{\mathcal{S}_{G}^{p}}:=\left(\widehat{\mathbb{E}}\Big[\sup_{t\in[0, T]}\vert\eta_{t}\vert^{p}\Big]\right)^{1/p}$. We denote by $\mathcal{S}_{G}^{p}\left(0, T\right)$ the completion of $\mathcal{S}_{G}^{0}\left(0, T\right)$ under the norm $\left\Vert\cdot\right\Vert_{\mathcal{S}_{G}^{p}}$.
\begin{Def}
For $\eta\in M_{G}^{0}\left(0, T\right)$ of the form \eqref{simp}, the It\^{o} integral with respect to $G$-Brownian motion is defined by the linear mapping $\mathcal{I}: M_{G}^{0}(0, T)\longrightarrow L_{G}^{2}(\Omega_{T})$,
\begin{equation*}
    \mathcal{I}(\eta):=\int_{0}^{T}\eta_{t}dB_{t}=\sum_{k=0}^{N-1}\xi_{k}(B_{t_{k+1}}-B_{t_{k}}),
\end{equation*}
which can be continuously extended to $\mathcal{I}: H_{G}^{1}(0,T)\longrightarrow L_{G}^{2}(\Omega_{T})$. On the other hand, the stochastic integral with respect to $(\langle B\rangle_t)_{t\geq 0}$ is defined by the linear mapping $\mathcal{Q}: M_{G}^{0}(0, T)\longrightarrow L_{G}^{1}(\Omega_{T})$,
\begin{equation*}
    \mathcal{Q}(\eta):=\int_{0}^{T}\eta_{t}d\langle B\rangle_{t}=\sum_{k=0}^{N-1}\xi_{k}(\langle B\rangle_{t_{k+1}}-\langle B\rangle_{t_{k}}),
\end{equation*}
which can be continuously extended to $\mathcal{Q}: H_{G}^{1}(0,T)\longrightarrow L_{G}^{1}(\Omega_{T})$.
\end{Def}
\begin{Lem}\label{BDG}\citep[BDG type inequality, see][Theorem~2.1]{Gao2009}.
Let $p\geq 2$, $\eta\in H_{G}^{p}(0, T)$ and $0\leq s\leq t\leq T$. Then,
\begin{eqnarray*}
  &&c_{p}\underline{\sigma}^{p}\widehat{\mathbb{E}}
    \left[\left(\int_{0}^{T}|\eta_{s}|^{2}ds\right)^{p/2}\right]  \\
  && \leq
    \widehat{\mathbb{E}}\left[\sup_{0\leq t\leq T}\left\vert\int_{0}^{t}
    \eta_{r}dB_{r}\right\vert^{p}\right]\leq
    C_{p}\overline{\sigma}^{p}\widehat{\mathbb{E}}
    \left[\left(\int_{0}^{T}|\eta_{s}|^{2}ds\right)^{p/2}\right],
\end{eqnarray*}
where $0<c_{p}<C_{p}<\infty $ are constants independent of $\eta$, $\underline{\sigma}$ and $\overline{\sigma}$.
\end{Lem}

For $\xi\in Lip(\Omega_{T})$, let
\begin{equation*}
    \mathcal{E}(\xi):=\widehat{\mathbb{E}}\left(\sup_{t\in[0, T]}\widehat{\mathbb{E}}_{t}[\xi]\right).
\end{equation*}
$\mathcal{E}$ is called the $G$-evaluation.

For $p\geq 1$ and $\xi\in Lip(\Omega_{T})$, define
\begin{equation*}
    \Vert\xi\Vert_{p, \mathcal{E}}:=\left(\mathcal{E}[\vert\xi\vert^{p}]\right)^{1/p}
\end{equation*}
and denote by $L_{\mathcal{E}}^{p}(\Omega_{T})$ the completion of $Lip(\Omega_{T})$ under the norm $\Vert\cdot\Vert_{p, \mathcal{E}}$.

The following estimate will be used in this paper.
\begin{Thm}\label{thm2}\citep[See][]{Song2011}.
For any $\alpha\geq 1$ and $\delta>0$, we have $L_{G}^{\alpha+\delta}(\Omega_{T})\subset L_{\mathcal{E}}^{\alpha}(\Omega_{T})$.
More precisely, for any $1<\gamma<\beta:=(\alpha+\delta)/\alpha$, $\gamma\leq 2$ and for all $\xi\in Lip(\Omega_{T})$, we have
\begin{equation*}\label{}
    \widehat{\mathbb{E}}\Big[\sup_{t\in[0, T]}\widehat{\mathbb{E}}_{t}[\vert\xi\vert^{\alpha}]\Big]\leq
    C\Big\{(\widehat{\mathbb{E}}[\vert\xi\vert^{\alpha+\delta}])^{\alpha/(\alpha+\delta)}+
    (\widehat{\mathbb{E}}[\vert\xi\vert^{\alpha+\delta}])^{1/\gamma}\Big\},
\end{equation*}
where
\begin{equation*}
    C=\frac{\gamma}{\gamma-1}(1+14\sum_{i=1}^{\infty}i^{-\beta/\gamma}).
\end{equation*}
\end{Thm}
\begin{Rem}\label{rm1}
By $\frac{\alpha}{\alpha+\delta}<\frac{1}{\gamma}<1$, we have
\begin{equation*}\label{}
    \widehat{\mathbb{E}}\Big[\sup_{t\in[0, T]}\widehat{\mathbb{E}}_{t}[\vert\xi\vert^{\alpha}]\Big]\leq
    2C\Big\{(\widehat{\mathbb{E}}[\vert\xi\vert^{\alpha+\delta}])^{\alpha/(\alpha+\delta)}+
    \widehat{\mathbb{E}}[\vert\xi\vert^{\alpha+\delta}]\Big\}.
\end{equation*}
Set
\begin{equation*}
    C_{1}=2\inf\Big\{\frac{\gamma}{\gamma-1}(1+14\sum_{i=1}^{\infty}i^{-\beta/\gamma}):1<\gamma<\beta, \gamma\leq 2\Big\},
\end{equation*}
then
\begin{equation}\label{eq8}
    \widehat{\mathbb{E}}\Big[\sup_{t\in[0, T]}\widehat{\mathbb{E}}_{t}[\vert\xi\vert^{\alpha}]\Big]\leq
    C_{1}\Big\{(\widehat{\mathbb{E}}[\vert\xi\vert^{\alpha+\delta}])^{\alpha/(\alpha+\delta)}+
    \widehat{\mathbb{E}}[\vert\xi\vert^{\alpha+\delta}]\Big\},
\end{equation}
where $C_{1}$ is a constant only depending on $\alpha$ and $\delta$.
\end{Rem}

\section{Large deviations for $G$-SDEs}
\label{Sec:G-SDEs}
In this section, we present the large deviations for $G$-SDEs obtained by \citet{Gao2010}. The authors use discrete time approximation to obtain their results.

First, we recall the following notations on large deviations under a sublinear expectation.

Let $(\chi, d)$ be a Polish space. Let $\left(U^{\varepsilon},\;\varepsilon>0\right)$ be a family of measurable maps from $\Omega$ into $(\chi, d)$ and let $\delta(\varepsilon)$, $\varepsilon>0$ be a positive function satisfying $\delta(\varepsilon)\rightarrow 0$ as $\varepsilon\rightarrow 0$.

A nonnegative function $\mathcal{I}$ on $\chi$ is called to be (good) \emph{rate function} if $\{x:\;\mathcal{I}(x)\leq\alpha\}$ (its level set) is (compact) closed for all $0\leq\alpha<\infty$.

$\left\{\widehat{C}(U^{\varepsilon}\in\cdot)\right\}_{\varepsilon>0}$ is said to satisfy large deviation principle with speed $\delta(\varepsilon)$ and with rate function $\mathcal{I}$ if for any measurable closed subset $\mathcal{F}\subset\chi$,
\begin{equation*}\label{UBLD}
    \limsup_{\varepsilon\rightarrow 0}\delta(\varepsilon)\log \widehat{C}\left(U^{\varepsilon}\in \mathcal{F}\right)\leq-\inf_{x\in \mathcal{F}}\mathcal{I}(x),
\end{equation*}
and for any measurable open subset $\mathcal{O}\subset\chi$,
\begin{equation*}\label{LBLD}
    \liminf_{\varepsilon\rightarrow 0}\delta(\varepsilon)\log \widehat{C}\left(U^{\varepsilon}\in \mathcal{O}\right)\geq-\inf_{x\in \mathcal{O}}\mathcal{I}(x).
\end{equation*}

In \citet{Gao2010}, for any $\varepsilon>0$, the authors considered the following random perturbation SDEs driven by $d$-dimensional $G$-Brownian motion $B$
\begin{equation*}
    X_{t}^{x,\varepsilon} =x+\int_{0}^{t}b^{\varepsilon}(X_{r}^{x,\varepsilon})dr+\varepsilon\int_{0}^{t}h^{\varepsilon}(X_{r}^{x,\varepsilon})d\langle B, B\rangle_{r/\varepsilon}+\varepsilon\int_{0}^{t}\sigma^{\varepsilon}(X_{r}^{x,\varepsilon})dB_{r/\varepsilon}
\end{equation*}
where $\langle B, B\rangle$ is treated as a $d\times d$-dimensional vector,
\begin{equation*}
    b^{\varepsilon}=(b_{1}^{\varepsilon}, \ldots, b_{n}^{\varepsilon})^{\tau}:\;\mathbb{R}^{n}\longrightarrow\mathbb{R}^{n},\;
    \sigma^{\varepsilon}=(\sigma_{i, j}^{\varepsilon}):\;\mathbb{R}^{n}\longrightarrow\mathbb{R}^{n\times d}
\end{equation*}
and $h^{\varepsilon}:\;\mathbb{R}^{n}\longrightarrow\mathbb{R}^{n\times d^{2}}$.

Consider the following conditions:
\begin{description}
  \item[(H1)] $b^{\varepsilon}$, $\sigma^{\varepsilon}$ and $h^{\varepsilon}$ are uniformly bounded;
  \item[(H2)] $b^{\varepsilon}$, $\sigma^{\varepsilon}$ and $h^{\varepsilon}$ are uniformly Lipschitz continuous;
  \item[(H3)] $b^{\varepsilon}$, $\sigma^{\varepsilon}$ and $h^{\varepsilon}$ converge uniformly to $b:=b^{0}$, $\sigma:=\sigma^{0}$ and $h:=h^{0}$ respectively.
\end{description}

Let $\mathcal{C}([0, T],\mathbb{R}^{n})$ be the space of $\mathbb{R}^{n}$-valued continuous functions $\varphi$ on $[0, T]$ and $\mathcal{C}_{0}([0, T],\mathbb{R}^{n})$ the space of $\mathbb{R}^{n}$-valued continuous functions $\widetilde{\varphi}$ on $[0, T]$ with $\widetilde{\varphi}_{0}=0$.

Define
\begin{align*}
\mathbb{H}^{d}:=&\Big\{\phi\in\mathcal{C}_{0}([0,T],\mathbb{R}^{d}): \phi\;\text{is absolutely continuous and} \\
&\quad\quad\quad\quad\Vert\phi\Vert_{\mathbb{H}}^{2}:=\int_{0}^{T}\vert\phi^{\prime}(r)\vert^{2}dr<+\infty\Big\}, \\
\mathbb{A}:=&\Big\{\eta=\int_{0}^{t}\eta^{\prime}(r)dr;\; \eta^{\prime}: [0,T]\longrightarrow\mathbb{R}^{d\times d}\; \text{Borel measurable and} \\
&\quad\quad\quad\quad\eta^{\prime}(t)\in\Sigma\;\textrm{ for all } t\in[0, T]\Big\}.
\end{align*}

We recall the following result of a joint large deviation principle for $G$-Brownian motion and its quadratic variation process.
\begin{Thm}\citep[See][p.~2225]{Gao2010}.
$\left\{\widehat{C}\left((\varepsilon B_{t/\varepsilon}, \varepsilon\langle B\rangle_{t/\varepsilon})\mid_{t\in[0, T]}\;\in\cdot\right)\right\}_{\varepsilon>0}$ satisfies large deviation principle with speed $\varepsilon$ and with rate function
\begin{equation*}
J(\phi, \eta)=
\begin{cases}
\frac{1}{2}\int_{0}^{T}\langle \phi^{\prime}(r), (\eta^{\prime}(r))^{-1}\phi^{\prime}(r)\rangle dr,  &  \text{if } (\phi, \eta)\in\mathbb{H}^{d}\times\mathbb{A}, \\
+\infty,  &  \text{otherwise}.
\end{cases}
\end{equation*}
\end{Thm}

For any $(\phi, \eta)\in\mathbb{H}^{d}\times\mathbb{A}$, let $\Psi(\phi,\eta)\in\mathcal{C}([0, T],\mathbb{R}^{n})$ be the unique solution of the following ordinary differential equation (ODE in short)
\begin{eqnarray*}
  \Psi(\phi,\eta)(t) &=& x+\int_{0}^{t}b(\Psi(\phi,\eta)(r))dr+\int_{0}^{t}\sigma(\Psi(\phi,\eta)(r))\phi^{\prime}(r)dr \\
  && +\int_{0}^{t}h(\Psi(\phi,\eta)(r))\eta^{\prime}(r)dr.
\end{eqnarray*}
\begin{Thm}\label{l31}\citep[See][p.~2233]{Gao2010}. Let $(H1)$, $(H2)$ and $(H3)$ hold. Then for any closed subset $\mathcal{F}$ and any open subset $\mathcal{O}$ in $\left(\mathcal{C}_{0}([0, T],\mathbb{R}^{d}), \Vert\cdot\Vert\right)\times\left(\mathcal{C}_{0}([0, T],\mathbb{R}^{d\times d}), \Vert\cdot\Vert\right)\times\left(\mathcal{C}_{0}([0, T],\mathbb{R}^{n}), \Vert\cdot\Vert\right)$,
\begin{equation*}
    \limsup_{\varepsilon\rightarrow 0}\varepsilon\log\widehat{C}\left((\varepsilon B_{t/\varepsilon}, \varepsilon\langle B\rangle_{t/\varepsilon}, X_{t}^{x, \varepsilon}-x)\mid_{t\in[0, T]}\;\in \mathcal{F}\right)\leq-\inf_{(\phi, \eta, \psi)\in \mathcal{F}}\widehat{I}(\phi, \eta, \psi),
\end{equation*}
and
\begin{equation*}
    \liminf_{\varepsilon\rightarrow 0}\varepsilon\log\widehat{C}\left((\varepsilon B_{t/\varepsilon}, \varepsilon\langle B\rangle_{t/\varepsilon}, X_{t}^{x, \varepsilon}-x)\mid_{t\in[0, T]}\;\in \mathcal{O}\right)\geq-\inf_{(\phi, \eta, \psi)\in \mathcal{O}}\widehat{I}(\phi, \eta, \psi),
\end{equation*}
where
\begin{equation*}
\widehat{I}(\phi, \eta, \psi)=
\begin{cases}
J(\phi, \eta),  &  \text{if } (\phi, \eta)\in\mathbb{H}^{d}\times\mathbb{A},\, x+\psi=\Psi(\phi,\eta) \\
+\infty,  &  \text{otherwise}.
\end{cases}
\end{equation*}
\end{Thm}
For $0\leq\alpha<1$ given and $n\geq1$, for each $\psi\in\mathcal{C}_{0}([0, T],\mathbb{R}^{n})$, set
\begin{equation*}
    \Vert\psi\Vert_{\alpha}:=\sup_{s,t\in[0, T]}\frac{\vert\psi(s)-\psi(t)\vert}{\vert s-t\vert^{\alpha}}
\end{equation*}
and
\begin{equation*}
    \mathcal{C}^{\alpha}_{0}([0, T],\mathbb{R}^{n}):=\Big\{\psi\in\mathcal{C}_{0}([0,T],\mathbb{R}^{n}): \lim_{\delta\rightarrow0}\sup_{\vert s-t\vert<\delta}\frac{\vert\psi(s)-\psi(t)\vert}{\vert s-t\vert^{\alpha}}=0,\Vert\psi\Vert_{\alpha}<\infty\Big\}.
\end{equation*}
\begin{Thm}\label{l311}\citep[See][p.~2227]{Gao2010}. Let $0\leq\alpha<1/2$ and let $(H1)$, $(H2)$ and $(H3)$ hold.
Then for any closed subset $\mathcal{F}$ and any open subset $\mathcal{O}$ in $\left(\mathcal{C}^{\alpha}_{0}([0, T],\mathbb{R}^{n}), \Vert\cdot\Vert_{\alpha}\right)$,
\begin{equation*}
    \limsup_{\varepsilon\rightarrow 0}\varepsilon\log\widehat{C}\left((X_{t}^{x, \varepsilon}-x)\mid_{t\in[0, T]}\;\in \mathcal{F}\right)\leq-\inf_{\psi\in \mathcal{F}}I(\psi),
\end{equation*}
and
\begin{equation*}
    \liminf_{\varepsilon\rightarrow 0}\varepsilon\log\widehat{C}\left((X_{t}^{x, \varepsilon}-x)\mid_{t\in[0, T]}\;\in \mathcal{O}\right)\geq-\inf_{\psi\in \mathcal{O}}I(\psi),
\end{equation*}
where
\begin{equation*}\label{}
    I(\psi)=\inf\Big\{J(\phi,\eta):\psi=\Psi(\phi,\eta)-x\Big\}.
\end{equation*}
\end{Thm}
We immediately have the following result which will be used in the following section.
\begin{Coro}\label{l31} Let $(H1)$, $(H2)$ and $(H3)$ hold. Then for any closed subset $\mathcal{F}$ and any open subset $\mathcal{O}$ in $\mathcal{C}_{0}([0, T],\mathbb{R}^{n})$,
\begin{equation*}
    \limsup_{\varepsilon\rightarrow 0}\varepsilon\log\widehat{C}\left((X_{t}^{x, \varepsilon}-x)\mid_{t\in[0, T]}\;\in \mathcal{F}\right)\leq-\inf_{\widetilde{\varphi}\in \mathcal{F}}\Lambda(\widetilde{\varphi}),
\end{equation*}
and
\begin{equation*}
    \liminf_{\varepsilon\rightarrow 0}\varepsilon\log\widehat{C}\left((X_{t}^{x, \varepsilon}-x)\mid_{t\in[0, T]}\;\in \mathcal{O}\right)\geq-\inf_{\widetilde{\varphi}\in \mathcal{O}}\Lambda(\widetilde{\varphi}),
\end{equation*}
where
\begin{equation*}\label{}
    \Lambda(\widetilde{\varphi})=\inf\Big\{J(\phi,\eta):x+\widetilde{\varphi}=
    \Psi(\phi,\eta)\Big\}.
\end{equation*}
\end{Coro}

In the following section, we consider the following $G$-SDE: for every $s\leq t\leq T$, $x\in\mathbb{R}^{n}$,
\begin{equation}\label{eq0002}
X_{t}^{s,x,\varepsilon} =x+\int_{s}^{t}b(X_{r}^{s,x,\varepsilon})dr+\varepsilon\int_{s}^{t}h(X_{r}^{s,x,\varepsilon})d\langle B, B\rangle_{r}+\varepsilon\int_{s}^{t}\sigma(X_{r}^{s,x,\varepsilon})dB_{r},
\end{equation}
where $b$, $\sigma$ and $h$ are bounded. In order to use the large deviation principle obtained by \citet{Gao2010}, we will transform the G-SDE \eqref{eq0002} in the following form:
\begin{equation*}
\widetilde{X}_{t}^{s,x,\varepsilon} =x+\int_{s}^{t}b^{\varepsilon}(\widetilde{X}_{r}^{s,x,\varepsilon})dr+\varepsilon\int_{s}^{t}h^{\varepsilon}(\widetilde{X}_{r}^{s,x,\varepsilon})d\langle \widetilde{B}, \widetilde{B}\rangle_{r/\varepsilon}+\varepsilon\int_{s}^{t}\sigma^{\varepsilon}(\widetilde{X}_{r}^{s,x,\varepsilon})d\widetilde{B}_{r/\varepsilon},
\end{equation*}
where $\widetilde{B}_{t}:=\frac{1}{\sqrt{\varepsilon}}B_{t\varepsilon}$, $b^{\varepsilon}:=b$, $h^{\varepsilon}:=\varepsilon h$ and $\sigma^{\varepsilon}:=\sqrt{\varepsilon}\sigma$.

\section{Large deviations for $G$-BSDEs}
\label{Sec:G-BSDEs}
\citet{Hu2014a} obtained the existence, uniqueness and a priori estimates of the following backward stochastic differential equation driven by $G$-Brownian motion
\begin{equation}\label{GBSDE}
    Y_{t}=\xi+\int_{t}^{T}f(r, Y_{r}, Z_{r})dr+\int_{t}^{T}g(r, Y_{r}, Z_{r})d\langle B\rangle_{r}-\int_{t}^{T}Z_{r}dB_{r}-(K_{T}-K_{t}),
\end{equation}
where $K$ is a decreasing $G$-martingale, under standard Lispchitz conditions on $f(r, y, z)$, $g(r, y, z)$ in $(y, z)$ and the integrability condition on $\xi$. The unique solution of the BSDE~\eqref{GBSDE} is the triple $(Y, Z, K)$. The solution of an SDE is one process, say $X$. The solution of a "traditional" BSDE is a pair $(Y, Z)$, the solution of a BSDE driven by a $G$-Brownian motion is a triplet.

To establish large deviation principle for $G$-BSDEs, we consider the following forward-backward stochastic differential
equation driven by $G$-Brownian motion (we use \emph{Einstein convention}): for every $s\leq t\leq T$, $x\in\mathbb{R}^{n}$,
\begin{eqnarray}
d X_{t}^{s,x,\varepsilon}&=&b(X_{t}^{s,x,\varepsilon})dt+\varepsilon h_{ij}(X_{t}^{s,x,\varepsilon})d\langle B^{i}, B^{j}\rangle_{t}+\varepsilon\sigma_{j}(X_{t}^{s,x,\varepsilon})dB^{j}_{t}, X_{s}^{s,x,\varepsilon}=x,\nonumber \\
Y_{t}^{s,x,\varepsilon}&=&\Phi(X_{T}^{s,x,\varepsilon})+\int_{t}^{T}f(r,X_{r}^{s,x,\varepsilon},Y_{r}^{s,x,\varepsilon},Z_{r}^{s,x,\varepsilon})dr
\nonumber \\
&&+\int_{t}^{T}g_{ij}(r,X_{r}^{s,x,\varepsilon},Y_{r}^{s,x,\varepsilon},Z_{r}^{s,x,\varepsilon})d\langle B^{i}, B^{j}\rangle_{r}-\int_{t}^{T}Z_{r}^{s,x,\varepsilon}dB_{r}\nonumber \\
&&\quad\quad-(K_{T}^{s,x,\varepsilon}-K_{t}^{s,x,\varepsilon})\label{eq4},
\end{eqnarray}
where
\begin{equation*}
    b, h_{ij}, \sigma_{j}:\mathbb{R}^{n}\longrightarrow\mathbb{R}^{n};\;
    \Phi:\mathbb{R}^{n}\longrightarrow\mathbb{R};\;f, g_{ij}:[0, T]\times\mathbb{R}^{n}\times\mathbb{R}\times\mathbb{R}^{d}\longrightarrow\mathbb{R}
\end{equation*}
are deterministic functions and satisfy the following assumptions:
\begin{description}
  \item[(A0)] $b$, $\sigma$ and $h$ are bounded, i.e., there exists a constant $L>0$ such that
  \begin{equation*}
    \sup_{x\in\mathbb{R}^{n}}\max\Big\{\vert b(x)\vert, \Vert\sigma(x)\Vert_{HS}, \Vert h(x)\Vert_{HS}\Big\}\leq L,
  \end{equation*}
  where $\Vert A\Vert_{HS}:=\sqrt{\sum_{ij} a_{ij}^{2}}$ is the \emph{Hilbert-Schmidt} norm of a matrix $A=(a_{ij})$.
  \item[(A1)] $h_{ij}=h_{ji}$ and $g_{ij}=g_{ji}$ for $1\leq i, j\leq d$;
  \item[(A2)] $f$ and $g_{ij}$ are continuous in $t$;
  \item[(A3)] There exist a positive integer $m$ and a constant $L>0$ such that
  \begin{eqnarray*}
     && \vert b(x)-b(x')\vert + \sum_{i, j=1}^{d}\vert h_{ij}(x)-h_{ij}(x')\vert \\
     && \qquad\qquad + \sum_{j=1}^{d}\vert \sigma_{j}(x)-\sigma_{j}(x')\vert\leq L\vert x-x'\vert, \\
     && \vert \Phi(x)-\Phi(x')\vert\leq L(1+\vert x\vert^{m}+\vert x'\vert^{m})\vert x-x'\vert, \\
     && \vert f(t, x, y, z)-f(t, x', y', z')\vert + \sum_{i, j=1}^{d}\vert g_{ij}(t, x, y, z)-g_{ij}(t, x', y', z')\vert \\
     && \qquad\qquad \leq L\Big[(1+\vert x\vert^{m}+\vert x'\vert^{m})\vert x-x'\vert+\vert y-y'\vert+\vert z-z'\vert\Big].
  \end{eqnarray*}
\end{description}
It follows from \citet{Peng2010,Hu2014a} that, under the assumptions $\bf{(A0)-(A3)}$, the $G$-BSDE~\eqref{eq4} has a unique solution $\{(Y_{t}^{s,x,\varepsilon},Z_{t}^{s,x,\varepsilon},K_{t}^{s,x,\varepsilon}):
s\leq t\leq T\}$. Moreover, for any $\alpha>1$, we have $Y^{s,x,\varepsilon}\in\mathcal{S}^{\alpha}_{G}(0, T)$, $Z^{s,x,\varepsilon}\in H^{\alpha}_{G}(0, T)$ and $K^{s,x,\varepsilon}$ is a decreasing $G$-martingale with $K_{s}^{s,x,\varepsilon}=0$ and $K_{T}^{s,x,\varepsilon}\in L^{\alpha}_{G}(\Omega_{T})$.

We consider the following deterministic system: for every $s\leq t\leq T$, $x\in\mathbb{R}^{n}$,
\begin{eqnarray}
d\varphi_{t}^{s,x}&=&b(\varphi_{t}^{s,x})dt,\;\varphi_{s}^{s,x}=x,\nonumber \\
\psi_{t}^{s,x}&=&\Phi(\varphi_{T}^{s,x})+\int_{t}^{T}f(r,\varphi_{r}^{s,x},\psi_{r}^{s,x},0)dr\nonumber \\
&&\quad+2\int_{t}^{T}G(g(r,\varphi_{r}^{s,x},\psi_{r}^{s,x},0))dr\label{eq5}.
\end{eqnarray}

\begin{Lem}\label{l5}
Let $\textbf{(A0)}$, $\textbf{(A1)}$ and $\textbf{(A3)}$ hold. Then
\begin{enumerate}
  \item Let $p\geq 2$. For any $\varepsilon\in(0,1]$, there exists a constant $C_{p}>0$, independent of $\varepsilon$, such that
\begin{equation}\label{eq100}
\widehat{\mathbb{E}}\Big(\sup_{s\leq t\leq T}\vert X_{t}^{s,x,\varepsilon}-\varphi_{t}^{s,x}\vert^{p}\Big)\leq C_{p}\varepsilon^{p}.
\end{equation}
  \item Moreover, $\left\{\widehat{C}\left((X_{t}^{s, x, \varepsilon}-x)\mid_{t\in[s, T]}\;\in\cdot\right)\right\}_{\varepsilon>0}$ satisfies a large deviation principle with speed $\varepsilon$ and with rate function
\begin{equation*}\label{}
    \Lambda(\widetilde{\varphi})=\inf\Big\{J(\phi,\eta):x+\widetilde{\varphi}=
    \widehat{\Psi}(\phi,\eta)\Big\},
\end{equation*}
where $\widehat{\Psi}(\phi,\eta)\in\mathcal{C}([s, T],\mathbb{R}^{n})$ be the unique solution of the following ODE
\begin{equation*}
    \widehat{\Psi}(\phi,\eta)(t) = x+\int_{s}^{t}b(\widehat{\Psi}(\phi,\eta)(r))dr.
\end{equation*}
\end{enumerate}
\end{Lem}
\begin{pf}
\textrm{ }\\
  $1.$ Let $u\in[s, T]$, we have
  \begin{eqnarray*}
    X_{u}^{s, x, \varepsilon}-\varphi_{u}^{s, x} &=& \int_{s}^{u}\left(b(X_{r}^{s,x,\varepsilon})-b(\varphi_{r}^{s, x})\right)dr+\varepsilon\int_{s}^{u}h(X_{r}^{s,x,\varepsilon})d\langle B, B\rangle_{r} \\
    && +\varepsilon\int_{s}^{u}\sigma(X_{r}^{s,x,\varepsilon})dB_{r}.
  \end{eqnarray*}
  Then, there exists a constant $C_{p}>0$,
  \begin{eqnarray*}
    \vert X_{u}^{s,x,\varepsilon}-\varphi_{u}^{s, x}\vert^{p} &\leq& C_{p}\Big\{\int_{s}^{u}\vert b(X_{r}^{s,x,\varepsilon})-b(\varphi_{r}^{s, x})\vert^{p}dr \\
    && +\varepsilon^{p}\overline{\sigma}^{p}\int_{s}^{u}\Vert h(X_{r}^{s,x,\varepsilon})\Vert^{p}dr \\
    &&+\varepsilon^{p}\Big\vert\int_{s}^{u}\sigma(X_{r}^{s,x,\varepsilon})dB_{r}\Big\vert^{p}\Big\} \\
    &\leq& C_{p}\Big\{\int_{s}^{u}\vert X_{r}^{s,x,\varepsilon}-\varphi_{r}^{s, x}\vert^{p}dr+\varepsilon^{p} \\
    && +\varepsilon^{p}\Big\vert\int_{s}^{u}\sigma(X_{r}^{s,x,\varepsilon})dB_{r}\Big\vert^{p}\Big\}. \\
  \end{eqnarray*}
  For $t\in[s, T]$,
  \begin{eqnarray*}
    \sup_{s\leq u\leq t}\vert X_{u}^{s,x,\varepsilon}-\varphi_{u}^{s, x}\vert^{p} &\leq& C_{p}\Big\{\sup_{s\leq u\leq t}\int_{s}^{u}\vert X_{r}^{s,x,\varepsilon}-\varphi_{r}^{s, x}\vert^{p}dr+\varepsilon^{p} \\
    && +\varepsilon^{p}\sup_{s\leq u\leq t}\Big\vert\int_{s}^{u}\sigma(X_{r}^{s,x,\varepsilon})dB_{r}\Big\vert^{p}\Big\} \\
    &\leq& C_{p}\Big\{\int_{s}^{t}\sup_{s\leq u\leq r}\vert X_{u}^{s,x,\varepsilon}-\varphi_{u}^{s, x}\vert^{p}dr+\varepsilon^{p} \\
    && +\varepsilon^{p}\sup_{s\leq u\leq t}\Big\vert\int_{s}^{u}\sigma(X_{r}^{s,x,\varepsilon})dB_{r}\Big\vert^{p}\Big\}. \\
  \end{eqnarray*}
  So taking the $G$-expectation, it follows from the BDG inequality that
  \begin{equation*}
    \widehat{\mathbb{E}}[\sup_{s\leq u\leq t}\vert X_{u}^{s,x,\varepsilon}-\varphi_{u}^{s, x}\vert^{p}]\leq C_{p}\varepsilon^{p}+C_{p}\int_{s}^{t}\widehat{\mathbb{E}}[\sup_{s\leq u\leq r}\vert X_{u}^{s,x,\varepsilon}-\varphi_{u}^{s, x}\vert^{p}]dr.
  \end{equation*}
  Therefore, by Gronwall's inequality,
  \begin{equation*}
    \widehat{\mathbb{E}}\Big(\sup_{s\leq u\leq T}\vert X_{u}^{s,x,\varepsilon}-\varphi_{u}^{s,x}\vert^{p}\Big)\leq C_{p}\varepsilon^{p}.
  \end{equation*}
  $2.$ Set $\widetilde{B}_{t}=\frac{1}{\sqrt{\varepsilon}}B_{t\varepsilon}$. Thanks to Remark~\ref{l6}, $\widetilde{B}$ is a $G$-Brownian motion. Then, we have $B_{t}=\sqrt{\varepsilon}\widetilde{B}_{t/\varepsilon}, \langle B, B\rangle_{t}=\varepsilon\langle \widetilde{B}, \widetilde{B}\rangle_{t/\varepsilon}$. Therefore, by the uniqueness of the solution of the $G$-SDEs, it is easy to check that $\{X_{t}^{s,x,\varepsilon}: s\leq t\leq T\}$ is the solution of the following $G$-SDE:
  \begin{equation*}
    \widetilde{X}_{t}^{s,x,\varepsilon} =x+\int_{s}^{t}b^{\varepsilon}(\widetilde{X}_{r}^{s,x,\varepsilon})dr+\varepsilon\int_{s}^{t}h^{\varepsilon}(\widetilde{X}_{r}^{s,x,\varepsilon})d\langle \widetilde{B}, \widetilde{B}\rangle_{r/\varepsilon}+\varepsilon\int_{s}^{t}\sigma^{\varepsilon}(\widetilde{X}_{r}^{s,x,\varepsilon})d\widetilde{B}_{r/\varepsilon},
  \end{equation*}
  where $b^{\varepsilon}$, $h^{\varepsilon}$ and $\sigma^{\varepsilon}$ have already been defined at the end of Section~\ref{Sec:G-SDEs}. Therefore, in view of assumption $\textbf{(A0)}$, the proof follows by virtue of Corollary~\ref{l31}. $\hfill\square$
\end{pf}
\begin{Propo}\label{propo1}
Let $p\geq 2$. For any $\varepsilon\in(0,1]$, we have
\begin{equation}\label{eq9}
    \widehat{\mathbb{E}}\Big[\sup_{s\leq t\leq T}\vert X_{t}^{s,x,\varepsilon}\vert^{p}\Big]\leq C(1+\vert x\vert^{p}),
\end{equation}
where the constant $C$ depends on $L$, $G$, $p$, $n$ and $T$.
\end{Propo}
\begin{pf}
By Proposition~4.1 in \citet{Hu2014b}, there exists a constant $C>0$ such that
\begin{equation*}
    \widehat{\mathbb{E}}_{s}\Big[\sup_{s\leq t\leq T}\vert X_{t}^{s, x, \varepsilon}-x\vert^{p}\Big]\leq C\Big(1+\vert x\vert^{p}\Big).
\end{equation*}
Then
\begin{equation*}
    \widehat{\mathbb{E}}\Big[\sup_{s\leq t\leq T}\vert X_{t}^{s, x, \varepsilon}-x\vert^{p}\Big]\leq C\Big(1+\vert x\vert^{p}\Big),
\end{equation*}
which implies the desired result. $\hfill\square$
\end{pf}
\begin{Thm}\label{thm1}
Let $\bf{(A0)-(A3)}$ hold. For any $\varepsilon\in(0,1]$, there exists a constant $C>0$, independent of $\varepsilon$, such that
\begin{equation*}
    \widehat{\mathbb{E}}\Big(\sup_{s\leq t\leq T}\vert Y_{t}^{s,x,\varepsilon}-\psi_{t}^{s,x}\vert^{2}\Big)\leq C\varepsilon^{2}.
\end{equation*}
\end{Thm}
\begin{pf}
We consider the following $G$-BSDE: for every $s\leq t\leq T$, $x\in\mathbb{R}^{n}$,
\begin{eqnarray}\label{eq7}
Y_{t}^{s,x}&=&\Phi(\varphi_{T}^{s,x})+\int_{t}^{T}f(r,\varphi_{r}^{s,x},Y_{r}^{s,x},Z_{r}^{s,x})dr
\nonumber \\
&&+\int_{t}^{T}g_{ij}(r,\varphi_{r}^{s,x},Y_{r}^{s,x},Z_{r}^{s,x})d\langle B^{i}, B^{j}\rangle_{r}-\int_{t}^{T}Z_{r}^{s,x}dB_{r}\nonumber \\
&&\quad\quad-(K_{T}^{s,x}-K_{t}^{s,x}).
\end{eqnarray}
Let $M^{s,x}$ be the following decreasing $G$-martingale:
\begin{equation*}
    M_{t}^{s,x}:=\int_{s}^{t}g_{ij}(r, \varphi_{r}^{s,x}, \psi_{r}^{s,x}, 0)d\langle B^{i}, B^{j}\rangle_{r}-2\int_{s}^{t}G\left(g(r, \varphi_{r}^{s,x}, \psi_{r}^{s,x}, 0)\right)dr.
\end{equation*}
Thanks to equation \eqref{eq5} and the uniqueness of the solution of the $G$-BSDEs, it is easy to check that $\{(\psi_{t}^{s,x}, 0, M_{t}^{s,x}): s\leq t\leq T\}$ is the solution of the $G$-BSDE \eqref{eq7}.

So, by Proposition~2.16 in \citet{Hu2014b}, there exists a constant $C>0$ such that
\begin{eqnarray*}
  \widehat{\mathbb{E}}\Big[\sup_{s\leq t\leq T}\vert Y_{t}^{s,x,\varepsilon}-\psi_{t}^{s,x}\vert^{2}\Big] &\leq& C\Big\{\widehat{\mathbb{E}}\Big[\sup_{t\in[s, T]}\widehat{\mathbb{E}}_{t}[\vert\Phi(X_{T}^{s,x,\varepsilon})-\Phi(\varphi_{T}^{s,x})\vert^{2}]\Big] \\ \\
  && +\Big(\widehat{\mathbb{E}}\Big[\sup_{t\in[s, T]}\widehat{\mathbb{E}}_{t}\Big[\Big(\int_{s}^{T}\widehat{h}_{r}dr\Big)^{4}\Big]\Big]\Big)^{1/2} \\
  && + \widehat{\mathbb{E}}\Big[\sup_{t\in[s, T]}\widehat{\mathbb{E}}_{t}\Big[\Big(\int_{s}^{T}\widehat{h}_{r}dr\Big)^{4}\Big]\Big]\Big\},
\end{eqnarray*}
where
\begin{eqnarray*}
    \widehat{h}_{r}&=&\vert f(r, X_{r}^{s, x, \varepsilon}, \psi_{r}^{s, x}, 0)-f(r, \varphi_{r}^{s, x}, \psi_{r}^{s, x}, 0)\vert \\
    &&\quad + \sum_{i, j=1}^{d}\vert g_{ij}(r, X_{r}^{s, x, \varepsilon}, \psi_{r}^{s, x}, 0)-g_{ij}(r, \varphi_{r}^{s, x}, \psi_{r}^{s, x}, 0)\vert.
\end{eqnarray*}
Therefore, in view of assumption $\textbf{(A3)}$, we have
\begin{eqnarray}\label{eq13}
  &&\widehat{\mathbb{E}}\Big[\sup_{s\leq t\leq T}\vert Y_{t}^{s,x,\varepsilon}-\psi_{t}^{s,x}\vert^{2}\Big]\nonumber \\
  &\leq& C\Big\{\widehat{\mathbb{E}}\Big[\sup_{t\in[s, T]}\widehat{\mathbb{E}}_{t}[(1+\vert X_{T}^{s,x,\varepsilon}\vert^{m}+\vert\varphi_{T}^{s,x}\vert^{m})^{2}\vert X_{T}^{s,x,\varepsilon}-\varphi_{T}^{s,x}\vert^{2}]\Big]\nonumber \\
  && +\Big(\widehat{\mathbb{E}}\Big[\sup_{t\in[s, T]}\widehat{\mathbb{E}}_{t}\Big[\Big(\int_{s}^{T}(1+\vert X_{r}^{s,x,\varepsilon}\vert^{m}+\vert\varphi_{r}^{s,x}\vert^{m})\vert X_{r}^{s,x,\varepsilon}-\varphi_{r}^{s,x}\vert dr\Big)^{4}\Big]\Big]\Big)^{1/2}\nonumber \\
  && + \widehat{\mathbb{E}}\Big[\sup_{t\in[s, T]}\widehat{\mathbb{E}}_{t}\Big[\Big(\int_{s}^{T}(1+\vert X_{r}^{s,x,\varepsilon}\vert^{m}+\vert\varphi_{r}^{s,x}\vert^{m})\vert X_{r}^{s,x,\varepsilon}-\varphi_{r}^{s,x}\vert dr\Big)^{4}\Big]\Big]\Big\},\nonumber \\
  &\leq& C\Big\{D_{1}+D_{2}^{\frac{1}{2}}+D_{2}\Big\},
\end{eqnarray}
where
\begin{eqnarray*}
  D_{1} &=& \widehat{\mathbb{E}}\Big[\sup_{t\in[s, T]}\widehat{\mathbb{E}}_{t}[(1+\vert X_{T}^{s,x,\varepsilon}\vert^{m}+\vert\varphi_{T}^{s,x}\vert^{m})^{2}\vert X_{T}^{s,x,\varepsilon}-\varphi_{T}^{s,x}\vert^{2}]\Big], \\
  D_{2} &=& \widehat{\mathbb{E}}\Big[\sup_{t\in[s, T]}\widehat{\mathbb{E}}_{t}\Big[\Big(\int_{s}^{T}(1+\vert X_{r}^{s,x,\varepsilon}\vert^{m}+\vert\varphi_{r}^{s,x}\vert^{m})\vert X_{r}^{s,x,\varepsilon}-\varphi_{r}^{s,x}\vert dr\Big)^{4}\Big]\Big].
\end{eqnarray*}
By Theorem~\ref{thm2} and \eqref{eq8} in Remark~\ref{rm1}, for any $\delta_{1}>0$, we get
\begin{eqnarray*}
  D_{1} &\leq& C_{1}\Big\{\Big(\widehat{\mathbb{E}}\Big[(1+\vert X_{T}^{s,x,\varepsilon}\vert^{m}+\vert\varphi_{T}^{s,x}\vert^{m})^{2+\delta_{1}}\vert X_{T}^{s,x,\varepsilon}-\varphi_{T}^{s,x}\vert^{2+\delta_{1}}\Big]\Big)^{\frac{2}{2+\delta_{1}}} \\
  && +\widehat{\mathbb{E}}\Big[(1+\vert X_{T}^{s,x,\varepsilon}\vert^{m}+\vert\varphi_{T}^{s,x}\vert^{m})^{2+\delta_{1}}\vert X_{T}^{s,x,\varepsilon}-\varphi_{T}^{s,x}\vert^{2+\delta_{1}}\Big]\Big\} \\
  &=& C_{1}\Big\{D_{1, 1}+D_{1, 2}\Big\}.
\end{eqnarray*}
Similarly, for any $\delta_{2}>0$, we get
\begin{eqnarray*}
  D_{2} &\leq& C_{2}\Big\{\Big(\widehat{\mathbb{E}}\Big[\Big(\int_{s}^{T}(1+\vert X_{r}^{s,x,\varepsilon}\vert^{m}+\vert\varphi_{r}^{s,x}\vert^{m})\vert X_{r}^{s,x,\varepsilon}-\varphi_{r}^{s,x}\vert dr\Big)^{4+\delta_{2}}\Big]\Big)^{\frac{4}{4+\delta_{2}}} \\
  &&+\widehat{\mathbb{E}}\Big[\Big(\int_{s}^{T}(1+\vert X_{r}^{s,x,\varepsilon}\vert^{m}+\vert\varphi_{r}^{s,x}\vert^{m})\vert X_{r}^{s,x,\varepsilon}-\varphi_{r}^{s,x}\vert dr\Big)^{4+\delta_{2}}\Big]\Big\} \\
  &=& C_{2}\Big\{D_{2, 1}+D_{2, 2}\Big\}.
\end{eqnarray*}
Then, by H\"{o}lder's inequality, \eqref{eq100} in Lemma~\ref{l5} and \eqref{eq9} in Proposition~\ref{propo1}, we can get
\begin{eqnarray*}
  D_{1,2} &\leq& C\Big(\widehat{\mathbb{E}}\Big[(1+\vert X_{T}^{s,x,\varepsilon}\vert^{m}+\vert\varphi_{T}^{s,x}\vert^{m})^{4+2\delta_{1}}\Big]\Big)^{1/2} \\
  &&\quad\times\Big(\widehat{\mathbb{E}}\Big[\vert X_{T}^{s,x,\varepsilon}-\varphi_{T}^{s,x}\vert^{4+2\delta_{1}}\Big]\Big)^{1/2} \\
  &\leq & C\varepsilon^{2+\delta_{1}}.
\end{eqnarray*}
Thus
\begin{equation}\label{eq11}
    D_{1}\leq C^{1}(\varepsilon^{2}+\varepsilon^{2+\delta_{1}}).
\end{equation}
Furthermore
\begin{eqnarray*}
  &&\widehat{\mathbb{E}}\Big[\Big(\int_{s}^{T}(1+\vert X_{r}^{s,x,\varepsilon}\vert^{m}+\vert\varphi_{r}^{s,x}\vert^{m})\vert X_{r}^{s,x,\varepsilon}-\varphi_{r}^{s,x}\vert dr\Big)^{4+\delta_{2}}\Big] \\
  &\leq& C\widehat{\mathbb{E}}\Big[\Big(1+\sup_{r\in[s, T]}\vert X_{r}^{s,x,\varepsilon}\vert^{m}+\sup_{r\in[s, T]}\vert\varphi_{r}^{s,x}\vert^{m}\Big)^{4+\delta_{2}}\Big(\sup_{r\in[s, T]}\vert X_{r}^{s,x,\varepsilon}-\varphi_{r}^{s,x}\vert\Big)^{4+\delta_{2}}\Big] \\
  &\leq& C\Big\{\Big(\widehat{\mathbb{E}}\Big[\Big(1+\sup_{r\in[s, T]}\vert X_{r}^{s,x,\varepsilon}\vert^{m}+\sup_{r\in[s, T]}\vert\varphi_{r}^{s,x}\vert^{m}\Big)^{8+2\delta_{2}}\Big]\Big)^{1/2} \\
  &&\times\Big(\widehat{\mathbb{E}}\Big[\sup_{r\in[s, T]}\vert X_{r}^{s,x,\varepsilon}-\varphi_{r}^{s,x}\vert^{8+2\delta_{2}}\Big]\Big)^{1/2}\Big\}
\end{eqnarray*}
Therefore
\begin{equation}\label{eq12}
    D_{2}\leq C^{2}(\varepsilon^{4}+\varepsilon^{4+\delta_{2}}).
\end{equation}
So, by virtue of \eqref{eq13}, \eqref{eq11} and \eqref{eq12}, we have
\begin{equation*}
    \widehat{\mathbb{E}}\Big[\sup_{s\leq t\leq T}\vert Y_{t}^{s,x,\varepsilon}-\psi_{t}^{s,x}\vert^{2}\Big]
    \leq C\varepsilon^{2}\Big(1+\varepsilon^{\delta_{1}/2}+1+\varepsilon^{\delta_{2}/2}+\varepsilon^{2}+\varepsilon^{2+\delta_{2}}\Big),
\end{equation*}
which leads to the end of the proof. $\hfill\square$
\end{pf}
We have an immediate consequence of Theorem~\ref{thm1}.
\begin{Coro}\label{coro1}
For any $\varepsilon\in(0,1]$ and all $x$ in a compact subset of $\mathbb{R}^{n}$, there exists a constant $C>0$, independent of $s$, $x$ and $\varepsilon$, such that
\begin{equation*}
\widehat{\mathbb{E}}\Big(\sup_{s\leq t\leq T}\vert Y_{t}^{s,x,\varepsilon}-\psi_{t}^{s,x}\vert^{2}\Big)\leq C\varepsilon^{2}.
\end{equation*}
\end{Coro}
\begin{Thm}\label{thm3}
Let $\bf{(A0)-(A3)}$ hold. For any $\varepsilon\in(0,1]$, there exists a constant $C>0$, independent of $\varepsilon$, such that
\begin{equation*}
    \widehat{\mathbb{E}}\Big[\int_{s}^{T}\vert Z_{r}^{s,x,\varepsilon}\vert^{2}dr\Big]+\widehat{\mathbb{E}}\Big(\sup_{s\leq t\leq T}\vert K_{t}^{s,x,\varepsilon}-M_{t}^{s,x}\vert^{2}\Big)\leq C\varepsilon^{2},
\end{equation*}
where $M^{s,x}$ is the following decreasing $G$-martingale:
\begin{equation*}
    M_{t}^{s,x}=\int_{s}^{t}g_{ij}(r, \varphi_{r}^{s,x}, \psi_{r}^{s,x}, 0)d\langle B^{i}, B^{j}\rangle_{r}-2\int_{s}^{t}G\left(g(r, \varphi_{r}^{s,x}, \psi_{r}^{s,x}, 0)\right)dr.
\end{equation*}
\end{Thm}
\begin{pf}
Applying It\^{o}'s formula to $\vert Y_{t}^{s,x,\varepsilon}-\psi_{t}^{s,x}\vert^{2}$, we have
\begin{eqnarray*}
  &&\vert Y_{s}^{s,x,\varepsilon}-\psi_{s}^{s,x}\vert^{2}+\int_{s}^{T}\vert Z_{r}^{s,x,\varepsilon}\vert^{2}d\langle B\rangle_{r} \\
  &&=\vert\Phi(X_{T}^{s,x,\varepsilon})-\Phi(\psi_{T}^{s,x})\vert^{2} \\
  &&+2\int_{s}^{T}(Y_{r}^{s,x,\varepsilon}-\psi_{r}^{s,x})\Big(f(r,X_{r}^{s,x,\varepsilon},Y_{r}^{s,x,\varepsilon},Z_{r}^{s,x,\varepsilon})
  -f(r,\varphi_{r}^{s,x},\psi_{r}^{s,x},0)\Big)dr \\
  &&+2\int_{s}^{T}(Y_{r}^{s,x,\varepsilon}-\psi_{r}^{s,x})\Big(g(r,X_{r}^{s,x,\varepsilon},Y_{r}^{s,x,\varepsilon},Z_{r}^{s,x,\varepsilon})
  -g(r,\varphi_{r}^{s,x},\psi_{r}^{s,x},0)\Big)d\langle B\rangle_{r} \\
  &&-2\int_{s}^{T}(Y_{r}^{s,x,\varepsilon}-\psi_{r}^{s,x})Z_{r}^{s,x,\varepsilon}dB_{r} \\
  &&-2\int_{s}^{T}(Y_{r}^{s,x,\varepsilon}-\psi_{r}^{s,x})d(K_{r}^{s,x,\varepsilon}-M_{r}^{s,x}).
\end{eqnarray*}
Therefore, in view of assumption $\textbf{(A3)}$, we have
\begin{eqnarray*}
  &&\int_{s}^{T}\vert Z_{r}^{s,x,\varepsilon}\vert^{2}d\langle B\rangle_{r} \\
  &&\leq\vert\Phi(X_{T}^{s,x,\varepsilon})-\Phi(\psi_{T}^{s,x})\vert^{2} \\
  &&+2L(1+d^{2}\overline{\sigma}^{2})\int_{s}^{T}(1+\vert X_{r}^{s,x,\varepsilon}\vert^{m}+\vert\varphi_{r}^{s,x}\vert^{m})\vert Y_{r}^{s,x,\varepsilon}-\psi_{r}^{s,x}\vert\vert X_{r}^{s,x,\varepsilon}-\varphi_{r}^{s,x}\vert dr \\
  &&+2L(1+d^{2}\overline{\sigma}^{2})\int_{s}^{T}\vert Y_{r}^{s,x,\varepsilon}-\psi_{r}^{s,x}\vert^{2}dr \\
  &&+2L(1+d^{2}\overline{\sigma}^{2})\int_{s}^{T}\vert Y_{r}^{s,x,\varepsilon}-\psi_{r}^{s,x}\vert\vert Z_{r}^{s,x,\varepsilon}\vert dr \\
  &&+2\Big\vert\int_{s}^{T}(Y_{r}^{s,x,\varepsilon}-\psi_{r}^{s,x})Z_{r}^{s,x,\varepsilon}dB_{r}\Big\vert \\
  &&+2C\sup_{s\leq r\leq T}\vert Y_{r}^{s,x,\varepsilon}-\psi_{r}^{s,x}\vert\sup_{s\leq r\leq T}\vert K_{r}^{s,x,\varepsilon}-M_{r}^{s,x}\vert.
\end{eqnarray*}
On the other hand,
\begin{eqnarray*}
  (K_{t}^{s,x,\varepsilon}-M_{t}^{s,x}) &=& (Y_{t}^{s,x,\varepsilon}-\psi_{t}^{s,x})-(Y_{s}^{s,x,\varepsilon}-\psi_{s}^{s,x}) \\
  &&+\int_{s}^{t}\widehat{f}_{r}dr+\int_{s}^{t}\widehat{g}_{r}d\langle B\rangle_{r}-\int_{s}^{t}Z_{r}^{s,x,\varepsilon}dB_{r}.
\end{eqnarray*}
where
\begin{eqnarray*}
    \widehat{f}_{r}&=&\vert f(r, X_{r}^{s, x, \varepsilon}, \psi_{r}^{s, x}, 0)-f(r, \varphi_{r}^{s, x}, \psi_{r}^{s, x}, 0)\vert \\
    \widehat{g}_{r}&=&\sum_{i, j=1}^{d}\vert g_{ij}(r, X_{r}^{s, x, \varepsilon}, \psi_{r}^{s, x}, 0)-g_{ij}(r, \varphi_{r}^{s, x}, \psi_{r}^{s, x}, 0)\vert.
\end{eqnarray*}
Thus
\begin{eqnarray}\label{eq14}
  \vert K_{t}^{s,x,\varepsilon}-M_{t}^{s,x}\vert &\leq& \Big\{2\sup_{s\leq r\leq T}\vert Y_{r}^{s,x,\varepsilon}-\psi_{r}^{s,x}\vert+\int_{s}^{t}\widehat{f}_{r}dr\nonumber \\
  && +\int_{s}^{t}\widehat{g}_{r}d\langle B\rangle_{r}+\Big\vert\int_{s}^{t}Z_{r}^{s,x,\varepsilon}dB_{r}\Big\vert\Big\}.
\end{eqnarray}
Then
\begin{eqnarray*}
  &&\int_{s}^{T}\vert Z_{r}^{s,x,\varepsilon}\vert^{2}d\langle B\rangle_{r} \\
  &&\leq\quad C_{1}\sup_{s\leq r\leq T}\vert Y_{r}^{s,x,\varepsilon}-\psi_{r}^{s,x}\vert^{2} \\
  &&+\quad C_{2}\int_{s}^{T}(1+\vert X_{r}^{s,x,\varepsilon}\vert^{m}+\vert\varphi_{r}^{s,x}\vert^{m})\vert Y_{r}^{s,x,\varepsilon}-\psi_{r}^{s,x}\vert\vert X_{r}^{s,x,\varepsilon}-\varphi_{r}^{s,x}\vert dr \\
  &&+\quad C_{3}\int_{s}^{T}\vert Y_{r}^{s,x,\varepsilon}-\psi_{r}^{s,x}\vert\vert Z_{r}^{s,x,\varepsilon}\vert dr \\
  &&+\quad 2\Big\vert\int_{s}^{T}(Y_{r}^{s,x,\varepsilon}-\psi_{r}^{s,x})Z_{r}^{s,x,\varepsilon}dB_{r}\Big\vert \\
  &&+\quad 2TC_{4}\sup_{s\leq r\leq T}\vert Y_{r}^{s,x,\varepsilon}-\psi_{r}^{s,x}\vert\Big[\int_{s}^{T}\vert Z_{r}^{s,x,\varepsilon}\vert dr\Big] \\
  &&+\quad 2C\sup_{s\leq r\leq T}\vert Y_{r}^{s,x,\varepsilon}-\psi_{r}^{s,x}\vert\Big[\sup_{s\leq t\leq T}\Big\vert\int_{s}^{t}Z_{r}^{s,x,\varepsilon}dB_{r}\Big\vert\Big].
\end{eqnarray*}
Now, by Theorem~\ref{thm1}, the BDG inequality and Young's inequality, $2uv\leq \lambda u^{2}+\frac{v^{2}}{\lambda}$ for $\lambda>0$, we obtain that
\begin{equation*}
    \widehat{\mathbb{E}}\Big[\int_{s}^{T}\vert Z_{r}^{s,x,\varepsilon}\vert^{2}dr\Big]\leq C\varepsilon^{2}+\frac{4\widetilde{C}}{\lambda}\widehat{\mathbb{E}}\Big[\int_{s}^{T}\vert Z_{r}^{s,x,\varepsilon}\vert^{2}dr\Big].
\end{equation*}
Then, taking $\lambda$ such that $\lambda>4\widetilde{C}$, we deduce that
\begin{equation*}
    \widehat{\mathbb{E}}\Big[\int_{s}^{T}\vert Z_{r}^{s,x,\varepsilon}\vert^{2}dr\Big]\leq C\varepsilon^{2}.
\end{equation*}
From \eqref{eq14}, we have
\begin{eqnarray*}
  \sup_{s\leq t\leq T}\vert K_{t}^{s,x,\varepsilon}-M_{t}^{s,x}\vert^{2} &\leq& C\Big\{\sup_{s\leq r\leq T}\vert Y_{r}^{s,x,\varepsilon}-\psi_{r}^{s,x}\vert^{2}+\int_{s}^{T}\widehat{f}_{r}^{2}dr \\
  && +\int_{s}^{T}\widehat{g}_{r}^{2}dr+\sup_{s\leq t\leq T}\Big\vert\int_{s}^{t}Z_{r}^{s,x,\varepsilon}dB_{r}\Big\vert^{2}\Big\}.
\end{eqnarray*}
Therefore, by the same arguments as above, we get
\begin{equation*}
    \widehat{\mathbb{E}}\Big[\sup_{s\leq t\leq T}\vert K_{t}^{s,x,\varepsilon}-M_{t}^{s,x}\vert^{2}\Big]\leq C\varepsilon^{2}.
\end{equation*}
The proof is complete. $\hfill\square$
\end{pf}
\begin{Rem}
As a consequence of Theorems~\ref{thm1} and \ref{thm3}, we get
\begin{equation*}
    \widehat{\mathbb{E}}\Big[\sup_{s\leq t\leq T}\vert Y_{t}^{s,x,\varepsilon}-\psi_{t}^{s,x}\vert^{2}+\int_{s}^{T}\vert Z_{r}^{s,x,\varepsilon}\vert^{2}dr+\sup_{s\leq t\leq T}\vert K_{t}^{s,x,\varepsilon}-M_{t}^{s,x}\vert^{2}\Big]\leq C\varepsilon^{2},
\end{equation*}
where $C$ is a positive constant and then the solution $\{(Y_{t}^{s,x,\varepsilon},Z_{t}^{s,x,\varepsilon},K_{t}^{s,x,\varepsilon}):
s\leq t\leq T\}$ of the $G$-BSDE~\eqref{eq4} converges to $\{(\psi_{t}^{s,x},0,M_{t}^{s,x}):
s\leq t\leq T\}$ where $\psi^{s, x}$ is the solution of the following backward ODE:
\begin{equation*}
    \psi_{t}^{s,x}=\Phi(\varphi_{T}^{s,x})+\int_{t}^{T}f(r,\varphi_{r}^{s,x},\psi_{r}^{s,x},0)dr
    +2\int_{t}^{T}G(g(r,\varphi_{r}^{s,x},\psi_{r}^{s,x},0))dr,
\end{equation*}
and $M^{s,x}$ is the following decreasing $G$-martingale:
\begin{equation*}
    M_{t}^{s,x}=\int_{s}^{t}g_{ij}(r, \varphi_{r}^{s,x}, \psi_{r}^{s,x}, 0)d\langle B^{i}, B^{j}\rangle_{r}-2\int_{s}^{t}G\left(g(r, \varphi_{r}^{s,x}, \psi_{r}^{s,x}, 0)\right)dr.
\end{equation*}
\end{Rem}

We recall a very important result in large deviation theory, used to transfer a LDP from one space to another.
\begin{Lem}\label{l41}(Contraction principle). Let $\{\mu_{\varepsilon}\}_{\varepsilon>0}$ be a family of probability measures that satisfies the large deviation principle with a good rate function $\Lambda$ on a Hausdorff topological space $\chi$, and for $\varepsilon\in(0,1]$, let $f_{\varepsilon}:\;\chi\longrightarrow\Upsilon$ be continuous functions, with $(\Upsilon, d)$ a metric space. Assume that there exists a measurable map $f:\;\chi\longrightarrow\Upsilon$ such that for any compact set $\mathcal{K}\subset\chi$,
\begin{equation}\label{eq10}
    \limsup_{\varepsilon\rightarrow 0}\sup_{x\in\mathcal{K}}d\left(f_{\varepsilon}(x),\;f(x)\right)=0.
\end{equation}
Suppose further that $\{\mu_{\varepsilon}\}_{\varepsilon>0}$ is exponentially tight.
Then the family of probability measures $\{\mu_{\varepsilon}\circ f_{\varepsilon}^{-1}\}_{\varepsilon>0}$ satisfies the LDP in $\Upsilon$ with the good rate function
\begin{equation*}
    \Pi(y)=\inf\Big\{\Lambda(x): x\in\chi, y=f(x)\Big\}.
\end{equation*}
\end{Lem}
\begin{pf}
First, observe that the condition \eqref{eq10} implies that for any compact set $\mathcal{K}\subset\chi$, the function $f$ is continuous on $\mathcal{K}\subset\chi$ (consequently that $f$ is continuous everywhere).

Since $\{\mu_{\varepsilon}\}_{\varepsilon>0}$ is exponentially tight, for every $\alpha<\infty$, there exists a compact set $\mathcal{K}_{\alpha}\subset\chi$ such that
\begin{equation*}
    \limsup_{\varepsilon\rightarrow 0}\varepsilon\log\mu_{\varepsilon}(\mathcal{K}_{\alpha}^{c})<-\alpha.
\end{equation*}
For every $\delta>0$, set
\begin{equation*}
    \Gamma_{\varepsilon,\delta}=\{x\in\chi: d(f_{\varepsilon}(x), f(x))>\delta\}.
\end{equation*}
We have
\begin{equation*}
    \mu_{\varepsilon}(\Gamma_{\varepsilon,\delta})
    \leq\mu_{\varepsilon}(\Gamma_{\varepsilon,\delta}\cap\mathcal{K}_{\alpha})
    +\mu_{\varepsilon}(\mathcal{K}_{\alpha}^{c}).
\end{equation*}
Given $\delta>0$, the first term on the right is zero for $\varepsilon$ small enough, so that
\begin{equation*}
    \limsup_{\varepsilon\rightarrow 0}\varepsilon\log\mu_{\varepsilon}(\Gamma_{\varepsilon,\delta})
    \leq\limsup_{\varepsilon\rightarrow 0}\varepsilon\log\mu_{\varepsilon}(\mathcal{K}_{\alpha}^{c})<-\alpha
\end{equation*}
and letting $\alpha\rightarrow\infty$, we obtain
\begin{equation*}
    \limsup_{\varepsilon\rightarrow 0}\varepsilon\log\mu_{\varepsilon}(\Gamma_{\varepsilon,\delta})=-\infty.
\end{equation*}
Therefore, the lemma follows from Corollary~4.2.21 p.~133 in \citet{Dembo1998}. $\hfill\square$
\end{pf}

Now consider
\begin{equation}\label{eq1}
    u^{\varepsilon}(t,x)=Y_{t}^{t,x,\varepsilon},\; (t,x)\in[0, T]\times\mathbb{R}^{n}.
\end{equation}

In \citet{Hu2014b} it is shown that $u^{\varepsilon}$ is a viscosity solution of the following nonlinear partial differential equation (PDE in short):
\begin{equation*}\label{}
\begin{cases}
&\partial_{t}u^{\varepsilon}+\mathcal{L}^{\varepsilon}\left(D^{2}_{x}u^{\varepsilon}, D_{x}u^{\varepsilon}, u^{\varepsilon}, x, t\right)=0, \\
&u^{\varepsilon}(T,x)=\Phi(x),
\end{cases}
\end{equation*}
where
\begin{eqnarray*}
  \mathcal{L}^{\varepsilon}\left(D^{2}_{x}u^{\varepsilon}, D_{x}u^{\varepsilon}, u^{\varepsilon}, x, t\right) &=& G\left(H\left(D^{2}_{x}u^{\varepsilon}, D_{x}u^{\varepsilon}, u^{\varepsilon}, x, t\right)\right)+\langle b(x), D_{x}u^{\varepsilon}\rangle \\
  &&+f\left(t, x, u^{\varepsilon}, \langle\varepsilon\sigma_{1}(x), D_{x}u^{\varepsilon}\rangle, \ldots, \langle\varepsilon\sigma_{d}(x), D_{x}u^{\varepsilon}\rangle\right),
\end{eqnarray*}
and
\begin{eqnarray*}
  H_{ij}\left(D^{2}_{x}u^{\varepsilon}, D_{x}u^{\varepsilon}, u^{\varepsilon}, x, t\right) &=& \langle D^{2}_{x}u^{\varepsilon}\varepsilon\sigma_{i}(x), \varepsilon\sigma_{j}(x)\rangle +2\langle D_{x}u^{\varepsilon}, \varepsilon h_{ij}(x)\rangle \\
   &&+2g_{ij}\left(t, x, u^{\varepsilon}, \langle\varepsilon\sigma_{1}(x), D_{x}u^{\varepsilon}\rangle, \ldots, \langle\varepsilon\sigma_{d}(x), D_{x}u^{\varepsilon}\rangle\right)
\end{eqnarray*}

We define the following
\begin{equation}\label{eq2}
    u^{0}(t,x)=\psi_{t}^{t,x},\; (t,x)\in[0, T]\times\mathbb{R}^{n}.
\end{equation}
\begin{Propo}\label{l42}
For any $\varepsilon >0$ and all $x\in\mathbb{R}^{n}$,
\begin{equation*}
Y_{t}^{s,x,\varepsilon}=u^{\varepsilon}(t,X_{t}^{s,x,\varepsilon}),\; \forall t\in [s,T].
\end{equation*}
\end{Propo}
\begin{pf}
Using the Markov property of the $G$-SDE and the uniqueness of the solution $Y^{s,x,\varepsilon}$ of the $G$-BSDE \eqref{eq4} to show that
\begin{equation*}
Y_{r}^{s,x,\varepsilon}=Y_{r}^{t,X_{t}^{s,x,\varepsilon},\varepsilon},
\; s\leq t\leq r\leq T.
\end{equation*}
Taking $r=t$, we deduce that $Y_{t}^{s,x,\varepsilon}=u^{\varepsilon}(t,X_{t}^{s,x,\varepsilon})$, which leads to the end of the proof. $\hfill\square$
\end{pf}

Let $\mathcal{C}_{0, s}([s, T],\mathbb{R}^{n})$ be the space of $\mathbb{R}^{n}$-valued continuous functions $\widetilde{\varphi}$ on $[s, T]$ with $\widetilde{\varphi}_{s}=0$.

Let $s\in[0, T]$ and $\varepsilon\geq 0$. We define the mapping $F^{\varepsilon}:\;\mathcal{C}_{0,s}([s, T],\mathbb{R}^{n})\longrightarrow\mathcal{C}([s, T],\mathbb{R}^{n})$ by
\begin{equation}\label{eq3}
    F^{\varepsilon}(\widetilde{\varphi})=[t\longmapsto u^{\varepsilon}(t, x+\widetilde{\varphi}_{t})],\; s\leq t\leq T,\;\widetilde{\varphi}\in\mathcal{C}_{0,s}([s, T],\mathbb{R}^{n}),
\end{equation}
where $u^{\varepsilon}$ is given by \eqref{eq1} and $u^{0}$ by \eqref{eq2}.

By virtue of \eqref{eq3} and Proposition~\ref{l42}, for any $\varepsilon >0$ and all $x\in\mathbb{R}^{n}$, we have $Y^{s,x,\varepsilon}=F^{\varepsilon}\left(X^{s,x,\varepsilon}-x\right)$.

We have the following result of large deviations
\begin{Thm} Let $\bf{(A0)-(A3)}$ hold. Then for any closed subset $\mathcal{F}$ and any open subset $\mathcal{O}$ in $\mathcal{C}([s, T], \mathbb{R}^{n})$,
\begin{equation*}
    \limsup_{\varepsilon\rightarrow 0}\varepsilon\log\widehat{C}\left(Y^{s, x, \varepsilon}\in\mathcal{F}\right)\leq-\inf_{\psi\in\mathcal{F}}\Pi(\psi),
\end{equation*}
and
\begin{equation*}
    \liminf_{\varepsilon\rightarrow 0}\varepsilon\log\widehat{C}\left(Y^{s, x, \varepsilon}\in\mathcal{O}\right)\geq-\inf_{\psi\in\mathcal{O}}\Pi(\psi),
\end{equation*}
where
\begin{equation*}
    \Pi(\psi)=\inf\Big\{\Lambda(\widetilde{\varphi}): \psi_{t}=F^{0}(\widetilde{\varphi})(t)=u^{0}(t,x+\widetilde{\varphi}_{t}),t\in[s, T],\widetilde{\varphi}\in\mathcal{C}_{0,s}([s, T], \mathbb{R}^{n})\Big\}.
\end{equation*}
\end{Thm}
\begin{pf}
Since the family $\left\{\widehat{C}\left((X_{t}^{s, x, \varepsilon}-x)\mid_{t\in[s, T]}\;\in\cdot\right)\right\}_{\varepsilon>0}$ is exponentially tight (see Lemma~3.4 p.~2235 in \citet{Gao2010}), by virtue of Lemma~\ref{l41} (contraction principle) and Lemma~\ref{l5}, we just need to prove that $F^{\varepsilon}$, $\varepsilon>0$ are continuous and $\{F^{\varepsilon}\}_{\varepsilon>0}$ converges uniformly to $F^{0}$ on every compact subset of $\mathcal{C}_{0,s}([s, T], \mathbb{R}^{n})$, as $\varepsilon\rightarrow 0$.

\emph{Continuity of $F^{\varepsilon}$}:

Let $\varepsilon> 0$ and $\widetilde{\varphi}\in\mathcal{C}_{0,s}([s,T],\mathbb{R}^{n})$. Let $(\widetilde{\varphi}^{m})_{m}$ be a sequence in $\mathcal{C}_{0,s}([s,T],\mathbb{R}^{n})$ which converges to $\widetilde{\varphi}$ under the uniform norm.

We set $\varphi^{m}=x+\widetilde{\varphi}^{m},\;\varphi=x+\widetilde{\varphi}$. So, $(\varphi^{m})_{m}$ is a sequence in $\mathcal{C}([s,T],\mathbb{R}^{n})$ which converges to $\varphi$ under the uniform norm. Fix $\zeta > 0$. Since $\Vert \varphi^{m}-\varphi\Vert_{\infty}\rightarrow 0$, there exists $M>0$ such that,
\begin{equation}\label{eq6}
\Vert \varphi^{m}\Vert_{\infty}\leq M, \; \Vert \varphi\Vert_{\infty}\leq M.
\end{equation}
Since $u^{\varepsilon}$ is a continuous function in $[0, T]\times\mathbb{R}^{n}$, it follows that $u^{\varepsilon}$ is uniformly continuous in $[s,T]\times B(0, M)$ where $B(0,M)$ is the closed ball centered at the origin with radius $M$ in $\mathbb{R}^{n}$.
Therefore, there exists $\eta>0$ such that for $r_{1}, r_{2}\in[s, T]$ and $z_{1}, z_{2}\in B(0,M)$, $\vert r_{1}-r_{2}\vert<\eta$ and $\vert z_{1}-z_{2}\vert<\eta$, we have
\begin{equation*}
\vert u^{\varepsilon}(r_{1}, z_{1})-u^{\varepsilon}(r_{2}, z_{2})\vert\leq\zeta.
\end{equation*}
Since there exists $m_{0}$ such that $\forall$ $m\geq m_{0}$, $\Vert \varphi^{m}-\varphi\Vert_{\infty}\leq\eta$, in view of \eqref{eq6}, for any $r\in[s,T]$ and for all $m\geq m_{0}$, we have
\begin{equation*}
\varphi^{m}_{r}, \varphi_{r}\in B(0,M) \quad \text{and} \quad \vert u^{\varepsilon}(r,\varphi^{m}_{r})-u^{\varepsilon}(r,\varphi_{r})\vert\leq\zeta.
\end{equation*}
Thus
\begin{equation*}
\vert u^{\varepsilon}(r,x+\widetilde{\varphi}^{m}_{r})-u^{\varepsilon}(r,x+\widetilde{\varphi}_{r})\vert\leq\zeta.
\end{equation*}
So we conclude that $F^{\varepsilon}(\widetilde{\varphi}^{m})\rightarrow F^{\varepsilon}(\widetilde{\varphi})$, which proves the continuity of $F^{\varepsilon}$ at $\widetilde{\varphi}$.

\emph{Uniform convergence of $F^{\varepsilon}$}:

Let $\mathcal{K}$ be a compact subset of $\mathcal{C}_{0,s}([s,T],\mathbb{R}^{n})$ and let
\begin{equation*}
\mathfrak{L}=\{\varphi_{r}: \widetilde{\varphi}\in\mathcal{K}, \varphi=x+\widetilde{\varphi}, r\in[s,T]\}.
\end{equation*}
Obviously, $\mathfrak{L}$ is a compact subset of $\mathbb{R}^{n}$. Thanks to Corollary~\ref{coro1}, there exists a positive constant $C$ such that
\begin{eqnarray*}
\sup_{\widetilde{\varphi}\in\mathcal{K}}\Vert F^{\varepsilon}(\widetilde{\varphi})-F^{0}(\widetilde{\varphi})\Vert_{\infty}^{2}&=&
\sup_{\widetilde{\varphi}\in\mathcal{K}}\sup_{r\in[s,T]}\vert u^{\varepsilon}(r,x+\widetilde{\varphi}_{r})-u^{0}(r,x+\widetilde{\varphi}_{r})\vert^{2} \\
&=&
\sup_{\widetilde{\varphi}\in\mathcal{K}}\sup_{r\in[s,T]}\vert u^{\varepsilon}(r,\varphi_{r})-u^{0}(r,\varphi_{r})\vert^{2} \\
&=&\sup_{\widetilde{\varphi}\in\mathcal{K}}\sup_{r\in[s,T]}\vert Y_{r}^{r,\varphi_{r},\varepsilon}-\psi_{r}^{r,\varphi_{r}}\vert^{2} \\
&\leq&\sup_{x\in\mathfrak{L}}\sup_{r\in[s,T]}\vert Y_{r}^{r,x,\varepsilon}-\psi_{r}^{r,x}\vert^{2} \\
&\leq& C\varepsilon^{2}.
\end{eqnarray*}
Therefore the uniform convergence of the mapping $F^{\varepsilon}$ towards $F^0$ follows. $\hfill\square$
\end{pf}


\bibliography{LDGBSDEs_references}

\end{document}